\newcommand{\N}{\mathbb{N}}               

\newcommand{\R}{\mathbb{R}}

\newcommand{\A}{\mathscr{A}}
\newcommand{\E}{\mathscr{E}}
\newcommand{\m}{\mathfrak{m}}

\newcommand{\Diff}{{\rm Diff}}   
\mathchardef\varepsilon="010F
\mathchardef\epsilon="0122
\mathchardef\vartheta="0112
\mathchardef\theta="0123
\mathchardef\varrho="011A
\mathchardef\rho="0125
\mathchardef\varphi="011E     
\mathchardef\phi="0127
\renewcommand \emptyset \varnothing
\documentclass[12pt]{article}            
\usepackage[all]{xy}
\usepackage[latin1]{inputenc}        
\usepackage[dvips]{graphics}
\usepackage[dvips]{graphicx}
\usepackage{amsfonts}
\usepackage{amssymb}
\usepackage{amsmath}
\usepackage{amstext}
\usepackage{amsbsy}
\usepackage{amsopn}
\usepackage{amsthm}
\usepackage{amscd}
\usepackage{amsxtra}
\usepackage{mathrsfs}
\usepackage{upref}
\author{{\sc Gianmarco Capitanio}}
\date{}                     
\title{Simple tangential families and perestroikas of their envelopes
} 

\begin{document}        

\theoremstyle{plain}
\newtheorem*{theorem*}{\bf Theorem}
\newtheorem{theorem}{\bf Theorem}
\newtheorem*{conjecture*}{\bf Conjecture}
\newtheorem{lemma}{\bf Lemma}
\newtheorem{proposition}{\bf Proposition}
\newtheorem*{proposition*}{\bf Proposition}
\newtheorem{corollary}{\bf Corollary}
\newtheorem*{corollary*}{\bf Corollary}
\theoremstyle{definition}
\newtheorem*{definition*}{\bf Definition}
\newtheorem*{definitions*}{\bf Definitions}

\newtheorem{example}{\bf Example}
\newtheorem*{example*}{\bf Example}
\theoremstyle{remark}
\newtheorem*{remark*}{\bf Remark}
\newtheorem{remark}{\bf Remark}
\newtheorem*{remarks*}{\bf Remarks}
\newtheorem{remarks}{\bf Remarks}

\maketitle
\begin{abstract}
Tangential families are $1$-parameter families of rays emanating 
tangentially from smooth curves.
We classify tangential family germs up to Left-Right equivalence: 
we prove that there are two infinite series  
and four sporadic simple singularities of tangential family germs 
(in addition to two stable singularities).  
We give their normal forms and miniversal tangential deformations
(i.e., deformations among tangential families), 
and we describe the corresponding envelope perestroikas of small codimension. 
We also discuss envelope singularities of non simple tangential families.
\end{abstract} 



\section{Introduction}\label{sct:1}

A tangential family is a system of rays emanating tangentially from a
smooth curve.
Tangential families naturally arise in the Geometry of Caustics (see
e.g. \cite{arnold2001}) and in Differential Geometry. 
For instance, every smooth curve in a Riemannian surface defines the
tangential family of its tangent geodesics. 
Tangential family theory is a generalization of Envelope theory. 
Indeed, every $1$-parameter family of plane curves is tangential, 
with respect to every generic point of any (geometric) branch of its envelope.

In \cite{stable} we classified stable singularities of 
tangential family germs (under small tangential deformations,
i.e. deformations among tangential families), and we proved that their
envelopes are smooth or have a second order self-tangency. 
In \cite{io2002} and \cite{io200?} we considered tangential 
families with singular support (the singularity being a 
semicubic cusp). 

In this paper we classify simple singularities of tangential family
germs.   
We prove that there are two infinite series and four sporadic simple
tangential family singularities (in addition to the two stable singularities).  
We give their normal forms and miniversal tangential deformations, and
we describe perestroikas of small codimension 
occurring to envelopes of simple tangential family germs 
under small tangential deformations. 
We also study envelopes of non simple tangential family
germs.

Our study of tangential families is related to Envelope Theory and Singularity
Theory, namely to the theories of Thom and Arnold on envelopes of
families of plane curves (see \cite{arnold1976a}, \cite{arnold1976b}
and \cite{Thom}), and to several branches of Projection Theory (see
\cite{ST2}), namely the classifications of simple projections of
surfaces and of projections of generic surfaces with boundary
(due to J.W. Bruce, P.J. Giblin and V.V. Goryunov, 
see \cite{goryunov}, \cite{BG} and \cite{goryunov2}), in our case the
boundary being fixed to be contained in the projection critical set. 

\section{Simple tangential family germs}\label{sct:2}

Unless otherwise specified, all the objects considered below 
are supposed to be of class $\mathscr C^\infty$; by plane curve we
mean an embedded $1$-submanifold of $\R^2$.     

Let us consider a map $f:\R^2\rightarrow \R^2$ of the plane $\R^2$, 
whose coordinates are denoted by $\xi$ and $t$.
If $\partial_t f$ vanishes nowhere, then $f_\xi:=f(\xi,\cdot)$
parameterizes an immersed smooth curve $\Gamma_\xi$.  
Hence, $f$ defines the $1$-parameter family of
curves $\{\Gamma_\xi:\xi\in\R\}$. 

\begin{definition*}
The family parameterized by $f$ is a {\it tangential family} whenever
(1) the curve parameterized by $f(\cdot,0)$, called the {\it support} of
the family, is smooth and (2) for every $\xi\in\R$, 
the family curve $\Gamma_\xi$ is tangent to the support at $f(\xi,0)$.
\end{definition*}

In other terms, a tangential family is a fibration, 
whose base is the support and whose fibers are the curves
$\Gamma_\xi$.   
A {\it $p$-parameter tangential deformation} of a tangential
family $f$ is a map $F:\R^2\times \R^p \rightarrow \R^2$, 
such that $F_\lambda:=F(\cdot;\lambda)$ is a
tangential family for every $\lambda\in\R^p$ and $F_0=f$.  
For instance, the tangent geodesics to a curve in a
Riemannian surface form a tangential family. 
A perturbation of the metric induces a tangential deformation on this
family.

The {\it graph} of a tangential family  $\{\Gamma_\xi:\xi\in\R\}$ is the surface 
$\Phi:=\cup_{\xi\in\R} \{(\xi,\Gamma_\xi)\} \subset \R\times \R^2$.
Let us denote by $\pi:\R\times \R^2\rightarrow \R^2$ 
the projection $(\xi,P)\mapsto P$. 
The {\it criminant set} of the tangential family is the critical set
of $\pi|_\Phi$.  
The {\it envelope} is the apparent contour of the graph
in the plane (i.e. the critical value set of $\pi|_\Phi$).   
By the very definition, the support of the family belongs to its envelope. 

Below we will study tangential family germs, so we will consider their
local parameterizations as map germs $(\R^2,0)\rightarrow (\R^2,0)$.  
Graphs of tangential family germs are smooth.

Denote by $\Diff(\R^2,0)$ the group of diffeomorphism
germs of the plane keeping fixed the origin and by 
$\m_{\xi,t}$ the ring of function germs $(\R^2,0)\rightarrow (\R,0)$ 
in the variables $\xi$ and $t$. 
Then the group $\A:= \Diff(\R^2,0)\times
\Diff(\R^2,0)$ acts on the space $(\m_{\xi,t})^2$ by 
$(\phi,\psi)\cdot f:= \psi\circ f\circ \phi^{-1}$, where
$(\phi,\psi)\in \A$ and $f\in (\m_{\xi,t})^2$.
Two elements of $(\m_{\xi,t})^2$ are said to be {\it Left-Right equivalent}, or
{\it $\A$-equivalent}, if they belong to the same $\A$-orbit.  
Thus, envelopes of  $\A$-equivalent tangential families are
diffeomorphic. 

The standard definition of versality of a deformation can be
translated in the setting of tangential deformations. 
A tangential deformation $F$ of a tangential family germ $f$ is said to be  
{\it versal} if any tangential deformation of $f$ is
$\A$-equivalent to a tangential deformation induced from $F$.  
Such a deformation is said to be {\it miniversal} 
whenever the base dimension $p$ is minimal. 
This minimal number $p$, depending only on the $\A$-orbit of the family, 
is called the {\it tangential codimension} of the singularity. 
A tangential family germ is {\it stable} if every $\A$-versal
tangential deformation is trivial; it is {\it simple} if by arbitrary sufficiently
small tangential deformations of it, we can obtain only a 
finite number of singularities.

Let us recall that a map germ 
$f:(\R^2,0)\rightarrow (\R^2,0)$ defines, by the formula 
$f^*g:=g\circ f$, a homomorphism from the ring $\mathscr E_{x,y}$ of 
the function germs in the target to the ring $\mathscr E_{\xi,t}$ of
the function germs in the source. 
Hence, we can consider every $\mathscr E_{\xi,t}$-module as an
$\mathscr E_{x,y}$-module via this homomorphism. 
The {\it tangent space} at $f$ is the  $\mathscr
E_{x,y}$-module defined by 
$$Tf:= \langle\partial_\xi f,
\partial_t f\rangle_{\E_{\xi,t}}  + f^*(\mathscr E_{x,y})\cdot 
\R^2  \ . $$ 
The {\it codimension} of $f$ is
the dimension of the real vector space $\E_{\xi,t}^2/Tf$. 

We state now our main result, proven in Section \ref{sct:6}.

\begin{theorem}\label{thm:1}
Every simple tangential family germ is $\A$-equivalent to a tangential
family germ $(\xi+t,\phi)$, of support $y=0$, 
where $\phi$ is one of the function germs listed in the table below,
in which $c$ and $\tau$ are the codimension
and the tangential codimension of the singularity.   
A miniversal tangential deformation of the normal form 
is obtained adding $\sum_i \lambda_i e_i$ to $\phi$ according to the table.  
\begin{center}
\begin{tabular}{|| c | c | c | c | c ||} 
\hline \hline
{\small Singularity} & $\phi$ &$c$& $\{e_i\}$ 
&$\tau$\\ \hline\hline
$\rm I$& $t^2$&$0$& --- &$0$\\ \hline 
${\rm II}$ & $t^2\xi$&$1$ & --- &$0$\\ \hline \hline 
$S_{1,n}$& $t^2(t+\xi)+t^4+t^{2n+3}$&$n+1$& $t^3, t^5, \dots, t^{2n+1}$  &$n$\\ \hline 
$T_n$& $t^3+t^2(t+\xi)^{n+1}$&
$2n+1$ &$t^2, t^2\xi, \dots, t^2\xi^{n-1}$ &$n$\\\hline
$S_{2,2}$& $t^2(t+\xi)+t^5+t^6$&$3$& $t^3,t^4$&$2$\\\hline
$S^\pm_{2,3}$& $t^2(t+\xi)+t^5\pm t^9$&$4$&$t^3,t^4,t^6$&$3$\\ \hline
$S_{2,4}$& $t^2(t+\xi)+t^5$&$5$&$t^3,t^4,t^6,t^9$ &$4$\\ \hline\hline
\end{tabular}
\end{center}
\end{theorem}

\begin{remarks*} 
(1) This table is a part of the more general classification of the simple
projections of surfaces in the plane, due to V. V. Goryunov
(see \cite{goryunov}). 

\noindent(2) In \cite{stable} we proved that singularities ${\rm I}$ and ${\rm II}$
are all the stable singularities; the corresponding envelopes are
respectively smooth and have an order $2$ self-tangency.   
\end{remarks*}

A curve with a singularity of order $b/a$ is {\it non symmetric} with
respect to a smooth curve tangent to it at the singular point if it is 
contained into one
of the two domains cut off near the singular point by the smooth curve.

\begin{corollary*} 
The envelope of every simple tangential family germ, having not a
singularity ${\rm I}$, has two tangent
branches, one of which is its support.   
More precisely:
\begin{itemize}
\item[(i)] if the germ has a singularity $S_{1,n}$, the envelope 
  second branch has a $(2n+3)/2$-cusp at the tangency point, this cusp
  being non symmetric to the support; 
\item[(ii)] if the germ has a singularity
  $S_{2,2}$, $S^\pm_{2,3}$ or $S_{2,4}$, then  the
  envelope second branch has a singularity of order $5/3$ at its
  tangency point with the support; 
\item[(iii)] if the germ has a singularity $T_n$, the envelope
  second branch is smooth and it has a tangency of order $3n+2$ with the
  support.
\end{itemize}\end{corollary*}

\begin{figure}[h]
  \centering
   \scalebox{.37}{\input{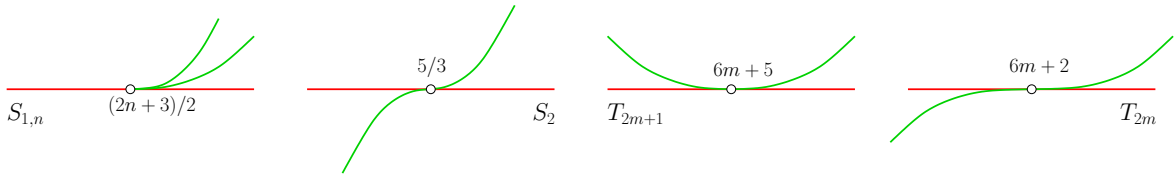}}
  \caption{Envelopes of simple tangential family germs.}
 \label{fig:1}
\end{figure}

Envelope singularities of simple tangential family germs 
are shown in Figure \ref{fig:1}.  

We end the classification of the simple singularities describing their adjacencies. 
A singularity $L$ is said to be {\it adjacent} to a singularity $K$, 
and we write $L\rightarrow K$, if every map germ in $L$ can be deformed
into a map germ in $K$ by an arbitrary small tangential deformation. 
If $L\rightarrow K\rightarrow K'$, then $L$ is also adjacent to 
$K'$; in this case we omit the arrow $L\rightarrow K'$. 
It follows from Theorem \ref{thm:1} that the adjacencies of the 
simple singularities of tangential family germs are as follows. 
$$\xymatrix{
{\rm I} & \ar[l] {\rm II} & \ar[l]T_{1}&\ar[l]T_{2}&\ar[l]T_{3}&\ar[l]T_{4}& 
          \ar[l]T_{5}&\ar[l] \dots \\
 &&\ar[lu]S_{1,1}&\ar[l]S_{1,2}&\ar[l]S_{1,3}&\ar[l]S_{1,4}& 
          \ar[l]S_{1,5}&\ar[l] \dots \\
&  &             &  \ar[lu]S_{2,2} & \ar[l]S^\pm_{2,3} & \ar[l]S_{2,4}\\  
}$$

\section{Bifurcation diagrams of simple singularities}\label{sct:3}

In this section we discuss the bifurcation diagrams of small
codimension simple singularities of tangential family germs, and the 
perestroikas of the corresponding envelopes. 

Let $L$ be a singularity, $f$ a map germ in $L$ and $F$ a miniversal 
tangential deformation of $f$. 
The {\it discriminant} of $L$ is, up to diffeomorphisms, the germ at
the origin of the set formed by 
the deformation parameter values $\lambda$ for which the envelope of $F_\lambda$ 
has more complicated singularities than for some arbitrary close value
of $\lambda$. 
The {\it bifurcation diagram} of a singularity is its discriminant, 
together with the envelopes of the deformed tangential family germs. 
In order to study the small codimension bifurctions diagrams of simple
singularities, we consider the miniversal tangential deformations
provided by Theorem \ref{thm:1}. 

We start with singularity $S_{1,n}$. 
By the corollary of Theorem \ref{thm:1}, the envelope
has two branches: the first is the family support $\gamma$, the second
has a $(2n+3)/2$-cusp tangent to $\gamma$. 

\begin{remark*} 
The discriminant of $S_{1,n}$ contains the flag 
$V_{n-1}\supset \cdots\supset V_0$, where 
$V_i$ is defined by $\lambda_1=...=\lambda_{n-i}=0$.  
A tangential family germ belonging to the stratum $V_i\smallsetminus
V_{i-1}$ has a singularity $S_{1,i}$, according to adjacency 
$S_{1,n}\rightarrow S_{1,n-1}$ (where $S_{1,0}:={\rm II}$). 
\end{remark*}

For singularity $S_{1,1}$, the envelope perestroika is shown in Figure \ref{fig:2}. 
\begin{figure}[h]     
\centering
  \scalebox{.3}{\input{figsimple2.pstex_t}}
  \caption{Envelope perestroika of singularity $S_{1,1}$.}
 \label{fig:2}
\end{figure} 

\begin{remark*}
The envelopes of the perturbed families form a singular surface germ
in $\R^3=\{x,y,\lambda\}$, diffeomorphic to a folded
umbrella with a cubically tangent smooth surface; the half line of the
umbrella self-intersections is tangent to the smooth surface. 
\end{remark*}

\begin{figure}[h]   
  \centering
  \scalebox{.33}{\input{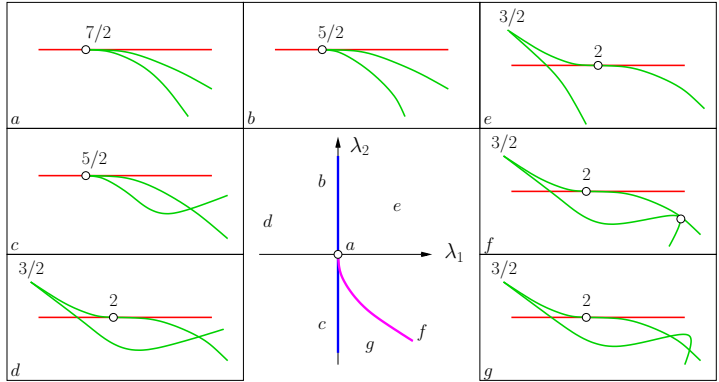}}
  \caption{Bifurcation diagram of singularity $S_{1,2}$.}
 \label{fig:4}
\end{figure}

The discriminant of singularity $S_{1,2}$ is represented in Figure
\ref{fig:4}.  
It has been found for me experimentally by Francesca Aicardi. 
Note that it contains the flag $V_1\supset V_0$ and a curve,
corresponding to a self-tangency of the envelope second branch.  

\begin{figure}[h]   
  \centering
  \scalebox{.34}{\input{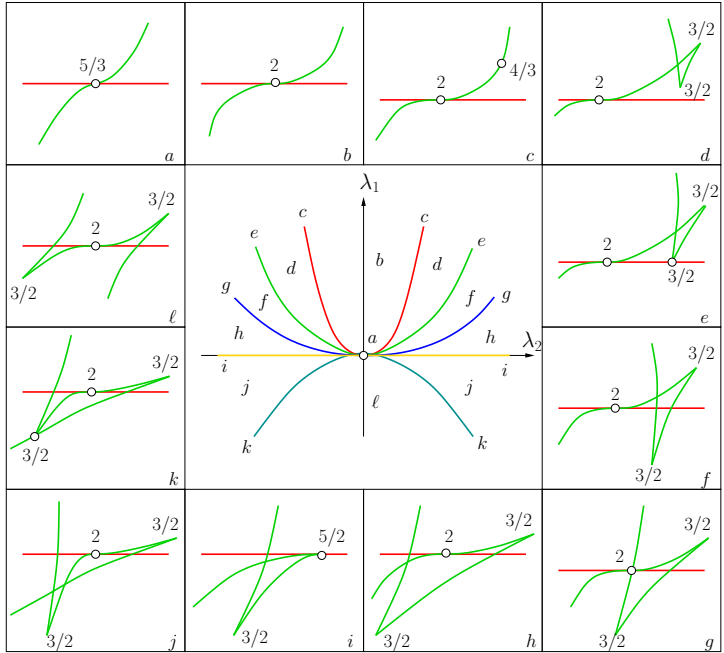}}
  \caption{Bifurcation diagram of the singularity $S_{2,2}$.}
 \label{fig:5}
\end{figure}

The $S_{2,2}$-discriminant in the plane $\{\lambda_1,\lambda_2\}$, 
shown in Figure \ref{fig:5}, is the union of the $\lambda_2$-axis and
four curve germs, tangent to it at the origin, the tangency order of 
all these curves being $1$.   
To my knowledge, this bifurcation diagram did not appear earlier in
Projection Theory. 

Finally, we consider singularities $T_n$. 

\begin{proposition} 
The discriminant of the singularity $T_n$ is an $n$-dimensional swallowtail.
\end{proposition}

\begin{proof}
The second branch of the envelope is the graph $y=\frac{4}{27} Q^3_\lambda(x)$, where 
$Q_\lambda(x):=x^{n+1}+\lambda_{n}x^{n-1}+ \dots +\lambda_1$,  
is a miniversal deformation of $x^{n+1}$ for the 
Left equivalence in the space of function
germs $(\R,0)\rightarrow\R$. 
The $T_n$-discriminant is formed by 
the values $\lambda$ for which $P_\lambda$ has roots of multiplicity
greater than $3$, i.e., for which $Q_\lambda$ has
multiple roots. 
Thus, it is a swallowtail.
\end{proof}

\begin{figure}[h]   
  \centering
  \scalebox{.4}{\input{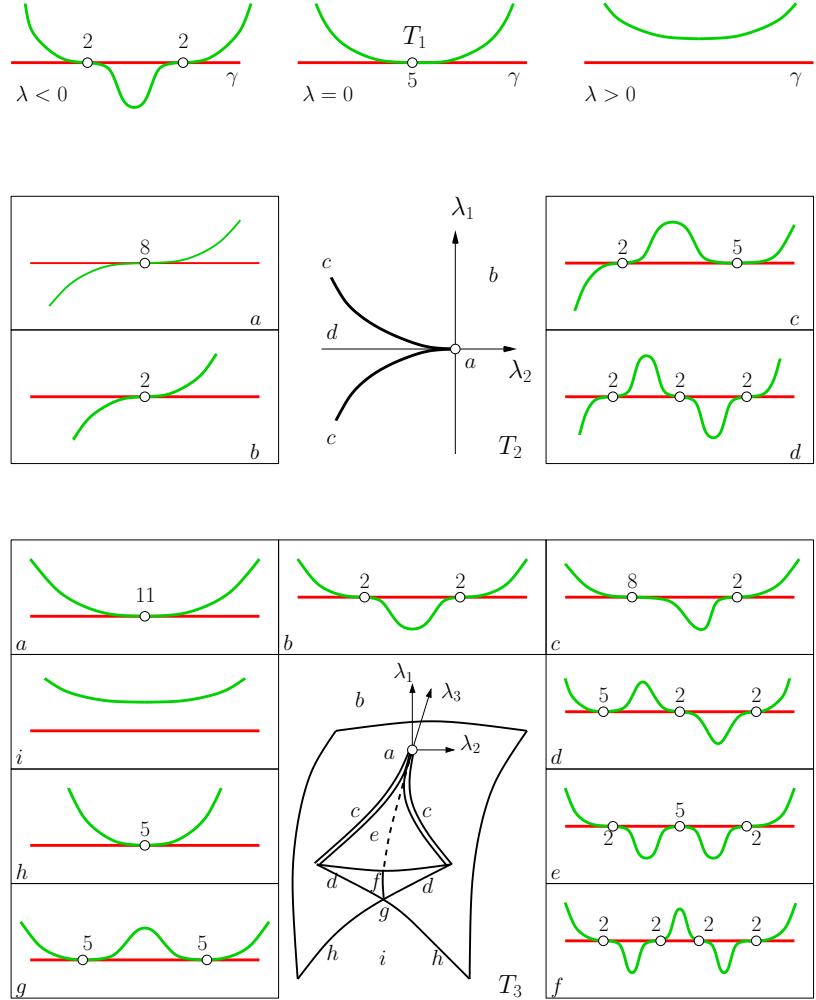}}
  \caption{Bifurcation diagrams of singularities $T_{n}$ for $n=1,2,3$.}
 \label{fig:6}
\end{figure}

Bifurcation diagrams of singularities $T_n$ are depicted in 
Figure \ref{fig:6} for $n=1,2,3$. 

\section{Classification of tangential family germs}\label{sct:4}

The set of tangential family germs is naturally decomposed according
to the configurtions of their criminant sets. 
We define now this decomposition.  

Given a tangential family germ, the fiber $\pi^{-1}(0,0)$ defines a 
{\it vertical direction} in the tangent plane to its graph at the origin.
The germ is said to be of {\it first type} if 
its criminant set has only one branch, of {\it second type} if it has exactly 
two branches and
these branches are smooth, non vertical and transversal each other.
In \cite{stable} we proved that a tangential family germ has a
singularity ${\rm I}$ (resp., ${\rm II}$) if and only if it is of
first type (resp., of second type).

\begin{definition*} 
Let $n\in\N\cup \{\infty\}$. 
A tangential family germ is said to be (1) {\it of type} $S_n$ 
if its criminant set has two branches, both smooth, one of which has
an order $n$ tangency with the vertical direction; (2) {\it of type}
$T_n$ if its  criminant set has two branches, both smooth, which have
an order $n$ tangency; (3) {\it of type} $U$ if its criminant set 
has more than two branches or at least a singular branch.
\end{definition*}
 
These singularity classes, in addition to stable singularities ${\rm I}$ and
${\rm II}$, cover all the tangential family singularities.    

\begin{definition*}
A tangential family germ of type $S_1$ is said to be {\it of type $S_{1,n}$} 
if its envelope has a cusp of order $(2n+3)/2$ at the origin.
\end{definition*}

\begin{theorem}\label{thm:4}
For every $n\in\N$, a tangential family germ has a singularity
$S_{1,n}$ (resp., $T_n$) if and only if it is of type $S_{1,n}$
(resp., $T_n$). 
A tangential family germ has a singularity $S_{2,2}$, $S^\pm_{2,3}$ or
$S_{2,4}$ if and only if it is of type $S_2$.
\end{theorem}

The next result describes envelopes of non simple tangential
family germs of type $S_{n\geq 3}$. 

\begin{theorem}\label{thm:5} 
For $n\in \N$, $n\geq 3$, the envelope of a $S_n$-type
tangential family germ has, in addition to the support, a branch 
with a singularity of order $(n+3)/(n+1)$, tangent to the support.  
If $n$ is even, this singularity is a cusp, non symmetric to the support.
\end{theorem}


Theorems \ref{thm:4} and \ref{thm:5} are proven in Sections
\ref{sct:5} and \ref{sct:6}.

We end this section with the description of the hierarchy of 
non simple singularities. 
We denote by $(L)$ and $(K)$ some classes of non simple
singularities. 
The arrow $(L)\rightarrow(K)$ means that there exists two
singularities $L'\subset(L)$ and $K'\subset(K)$ such that $L'\rightarrow K'$. 
The main adjacencies of non simple singularities are as follows.
$$
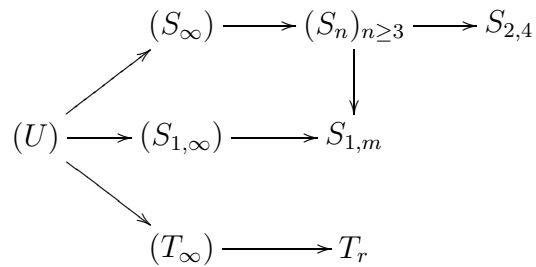
\xymatrix{
           &(S_{\infty}) \ar[r] & (S_n)_{n\geq 3}\ar[d] \ar[r]&S_{2,4}         \\
(U)\ar[r] \ar[ru] \ar[rd] &(S_{1,\infty}) \ar[r]& S_{1,m} &    \\
           &(T_\infty)    \ar[r]   & T_r& \\
}$$

\section{Preliminary results}\label{sct:5}

Here we introduce some technical results we will use in next section. 
We start describing standard parameterizations of tangential family
germs. 

Let us fix a tangential family germ.   
One easily verifies (see \cite{stable}) that in a coordinate
system in which the family support is the germ of $y=0$, 
the tangential family is parameterized by a map germ $(\xi+t,\phi)$, 
where $\phi$ is of the form $\alpha t^3 + t^2\sum_{i=0}^\infty
k_i \xi^i + t^3 \cdot \delta(0)+R(\xi)$, where $\alpha,k_i \in\R$, 
$\delta(n)$ denotes any function of two variables with vanishing
$n$-jet at the origin and $R$ is flat (i.e., its
Taylor expansion is zero).  
Such a parameterization is called {\it prenormal form}.
The prenormal form of a tangential family germ is not unique.

In  \cite{stable} we proved that {\it a tangential family germ 
is of first type if and only if $k_0\not=0$; it is of second type if and
only if $k_0=0$ and $k_1\not=0,\alpha$}.
Therefore, for any tangential family germ neither of first nor of
second type, we have $k_0=0$ and $k_1(k_1-\alpha)=0$.
 
\begin{lemma}\label{lemma:2}
Let $(\xi+t,\phi)$ be a prenormal form of a tangential family germ,
where $\phi$ is as above. 
The germ is of type $(S_n)_{n\in\N}$, $(T_n)_{n\in\N}$ 
or $U$ if and only $\alpha=k_1\not=0$, 
$\alpha\not=k_1=0$ or $\alpha=k_1=0$ respectively.  
More precisely, 
\begin{itemize}
\item[(i)] the family is of type $S_n$ if and only if 
$\phi(0,t)=\alpha t^3+O(t^{2n+3})$; 
\item[(ii)] the family is of type $T_n$ if and only if
  $k_0=\dots=k_n=0$ and $\alpha,k_{n+1}\not=0$,
\item[(iii)] the family is of type $U$ if and only if 
$\phi(\xi,t)=t^2\cdot\delta(1)$.
\end{itemize}
\end{lemma}

Theorem \ref{thm:5} follows from Lemma \ref{lemma:2} by explicit
computation of the envelope. 

\begin{remark}\label{rk:1}
Let $(\xi+t,\phi)$ the prenormal form of a tangential family germ. 
Then every $p$-parameter tangential deformation of it is equivalent to
a tangential deformation of the form $(\xi+t,\phi+t^2\psi)$, where
$\psi=\psi(\xi,t;\lambda)$ is a function germ
$(\R^2\times\R^p,0)\rightarrow (\R^2,0)$. 
\end{remark}

Let us consider now the space $\R[[\xi,t]]$, endowed with a 
quasihomogeneous filtration defined by weights $\deg(\xi)=a$
and $\deg(t)=b$, where $a$ and $b$ are coprime natural numbers. 
We denote by 
$\tilde\m^{p}_{\xi,t}$ the $\E_{\xi,t}$-ideal
generated by the monomials of weighted degree $p$, and by 
$\tilde \delta(p)$ any function germ with zero weighted $p$-jet. 

The {\it reduced tangent space} at $f$ is the $\mathscr
E_{x,y}$-module defined by 
$T_rf:= \mathfrak g_+ +\mathcal M$,  
where $\mathfrak g_+$ is the space of the vector field
germs having positive weighted order, 
$\mathcal M$ is the $\mathscr E_{x,y}$-module 
$f^*\left(\m_{x,y}^2 \oplus \R \cdot y\right) \times 
f^*\left(\m_{x,y}^2 \oplus \R \cdot x\right)$, 
and $\m^2_{x,y}$ is the second power of the maximal ideal
$\m_{x,y}$. 

\begin{remark}\label{rk:2}
Let $R$ be an element of $\tilde \m^p\times\tilde \m^q$, 
such that  $R\in T_rf/(\tilde \m^{p+1}\times\tilde \m^{q+1})$.
Then the weighted $(p,q)$-jet of $f$ is $\A$-equivalent to that of
$f+R$.    
\end{remark}

\section{Proofs}\label{sct:6}

Here we prove Theorems \ref{thm:1} 
and \ref{thm:4}.   
For the computations we skip, we refer to \cite{these}.
A {\it vector monomial} is by definition an element of
$(\m_{\xi,t})^2$ having a component which is a monomial, the other
component being zero. 

\bigskip\noindent{\bf $S_1$-type tangential families.}
Set $\deg(t)=1$, $\deg(\xi)=2$. 
Consider a tangential family germ of type $S_1$; by  Lemma
\ref{lemma:2}, its prenormal form is $\A$-equivalent to  
$(\xi,t^4+t^2\xi+t^2\cdot
 \tilde\delta(2))$.
It follows from the geometry of the Newton Diagram of 
$t^4+t^2\xi$ that the preceding map germ 
is formally $\A$-equivalent to 
$(\xi, t^4+t^2\xi)+\sum_{i=1}^\infty b_i (0, t^{2i+3})$, for some
$b_i\in\R$. 
Such a germ is of type $S_{1,n}$ if and only if 
$b_i=0$ for $i<n$ and $b_n\not=0$. 
In this case we may assume $b_n=1$ after rescaling. 

We set $f_n(\xi,t):=(\xi,t^4+t^2\xi+t^{2n+3})$. 
Since $T_rf_n$ contains 
$\tilde\m^{2n+2}_{\xi,t}\times \tilde\m^{2n+4}_{\xi,t}$ 
and the vector monomials $(0,t^{2\ell+4})$ for $\ell\in\N$, 
Remark \ref{rk:2} implies that every $S_{1,n}$ type tangential family germ is 
$\A$-equivalent to $f_n$ ($n\in\N$), and hence to the $S_{1,n}$-normal
form $f^\prime_n$ listed in Theorem \ref{thm:1}.  

The $\E_{\xi,t}$-ideal $\langle f_n\rangle$ is generated by $\xi$ and
$t^4$.
By the Preparation Theorem of Malgrange and Mather 
(see e.g. \cite{Martinet}), $\{1,t,t^2, t^3\}$ is a
generator system of the $\mathscr E_{x,y}$-module $\mathscr E_{\xi,t}$. 
Now, $Tf_n$ contains $\tilde\m^{2n}_{\xi,t}\times \tilde\m^{2n+2}_{\xi,t}$ and 
the vector monomials $(0,t^{2\ell})$ for $\ell\in\N$. 
Hence, the vectors $(0,t^{2i+1})$, $i=0,\dots,n$, form a basis of the 
real vector space $\E_{\xi,t}^2/Tf_n$. 
Therefore, the codimension of $f_n$ is $n+1$ and 
$f_n+\sum_{i=0}^n (0,\lambda_i t^{2i+1})$ 
is an $\A$-miniversal deformation of $f_n$. 
As a consequence, $f^\prime_n+\sum_{i=0}^n (0,\lambda_i t^{2i+1})$ is an
$\A$-miniversal deformation of the normal form $f^\prime_n$; moreover,
the singularities of class $S_{1,\infty}$ have infinite codimension. 

Let us consider a tangential deformation $F=F(\xi,t;\alpha)$ 
of the tangential family germ $f^\prime_n$. 
By the very definition of versality, there exist function germs
$\lambda_i(\alpha)$ such that $F$ is equivalent to 
$f^\prime_n+\sum_{i=0}^n (0,\lambda_i(\alpha) t^{2i+1})$. 
By Remark \ref{rk:1}, such a deformation is tangential if and only if
$\lambda_0\equiv 0$. 
Thus,  $f^\prime_n+\sum_{i=1}^n (0,\lambda_i t^{2i+1})$ is a
miniversal tangential deformation of $f^\prime_n$.

\bigskip\noindent{\bf $S_2$-type tangential families.}
Set $\deg(\xi)=3$, $\deg(t)=1$.
By Lemma \ref{lemma:2}, every tangential family germ of type $S_2$ is
$\A$-equivalent to a map germ $g+(0,\tilde\delta(5))$,
where $g(\xi,t):=(\xi,t^5+t^2\xi)$.
It follows from the inclusion 
$\tilde \m_{\xi,t}^8\times \tilde \m_{\xi,t}^{10}\subset T_rg$ 
and a well known Theorem of Gaffney (see \cite{Gaffney}, \cite{ST1})
that $g$ is $9$-$\A$-determined.
One checks that there are exactly four 
$\A$-orbits in $(\m_{\xi,t})^2$ over the weighted $5$-jet $g$, 
represented by the map germs $g_2:= g+(0,t^6)$, 
$g_3^\pm := g\pm(0,t^9)$ and $g_4:= g$ (the forthcoming computations
being identical for $g_3^+$ and for
$g_3^-$, we omit the sign $\pm$). 
These germs are $\A$-equivalent to the normal forms $g'_i$ given in
Theorem \ref{thm:1}.

In all the four cases we deal with, the quotient space $\E_{\xi,t}/\langle
g_i\rangle$, where $i=2,3,4$, is generated by
$\{t, t^2, t^3, t^4\}$ (Preparation Theorem). 
Set $e_0:=t$, $e_1:=t^3$, $e_2:=t^4$, $e_3:=t^6$, $e_4:=t^9$.
Using the inclusions $\tilde \m_{\xi,t}^3\times \tilde
\m_{\xi,t}^{5}\subset T_rg_2$, $\tilde \m_{\xi,t}^5\times \tilde
\m_{\xi,t}^{7}\subset T_rg_3$ and the above inclusion concerning
$g_4$, one proves
that $\{(0,e_j):j=1,\dots i\}$ is a basis of 
the real vector space $\E_{\xi,t}^2/Tg_i$, so 
the map germs  $g_i+ \sum_{j=0,...,i} \lambda_j (0,e_j)$ are 
$\A$-miniversal deformations of the map germs $g_i$; in
particular, the codimensions of $g_2$, $g_3$ and $g_4$ are
$3$, $4$ and $5$.
Thus, the map germs $g_i'+ \sum_{j=0,...,i} \lambda_j (0,e_j)$ are 
$\A$-miniversal deformations of the normal forms  $g'_i$
As for $S_{1,n}$-type tangential families, we see that these 
deformations are miniversal tangential deformations if and only if
$\lambda_0\equiv 0$. 

\bigskip
\noindent{\bf $T$-type tangential families.}
Set $\deg(t)=n+1$ and $\deg(\xi)=1$ and fix $n\in\N$.  
By Lemma \ref{lemma:2}, any $T_n$-type tangential
family germ is $\A$-equivalent to $h_n+(0,\tilde \delta(3n+3))$, 
where $h_n(\xi,t):=(\xi,t^3+t^2\xi^{n+1})$.
Since $\tilde\m_{\xi,t}^2\times \tilde\m_{\xi,t}^{3n+4} \subset
T_rh_n$, every such a germ is $\A$-equivalent to $h_n$, due to Remark
\ref{rk:2}, and then to the normal form $h^\prime_n$ listed in 
Theorem \ref{thm:1}.   

Now, $\langle h_n\rangle$ is generated by $\xi$ and $t^3$.   
Since $\E_{\xi,t}\times \tilde\m_{\xi,t}^{3n+2} \subset Th_n$,  
the Preparation Theorem implies that the vectors 
$(0,t\xi^{j}), j=0,\dots,n$ and $(0,t^2\xi^i), i=0,\dots,n-1$ form a
basis of the real vector space $\E_{\xi,t}^2/Tf_n$.  
Therefore, the codimension of $h_n$ is $2n+1$ (so, the codimension of
any $T_\infty$ singularity is infinite). 
Actually, we get still a basis if we replce the vectors $(0,t\xi^{j})$
and $(0,t^2\xi^i)$
by $v_j:= (0,t(\xi-t)^j)$ and $w_i:=j(0,t^2(\xi-t)^i)$, so 
$h_n+\sum_{i=1}^{n}\lambda_i w_i+\sum_{j=0}^n \mu_j v_j$ 
is an $\A$-miniversal deformation of $h_n$. 
Thus, 
$h^\prime_n+\sum_{i=1}^{n}\lambda_i (0,t^2\xi^{i-1})+\sum_{j=0}^n \mu_j
(0,t\xi^j)$
is an $\A$-miniversal deformation of the normal form $h^\prime_n$. 
Any tangential deformation of the $T_n$-normal form is equivalent to
the above map germ, in which $\lambda_j$ and $\mu_i$ are function
germs depending on the parameter of the tangential deformation. 
Due to Remark \ref{rk:1}, all the function germs $\mu_i$ have to be
identically zero. 

\bigskip
The above computations end the proof of Theorem \ref{thm:4}.
To complete the proof of Theorem \ref{thm:1}, it remains to show that the
singularity classes $S_{1,\infty}$, $T_{\infty}$, $S_{n\geq3}$ and $U$
are non simple.
For singularities of classes $S_{1,\infty}$ and
$T_\infty$, this follows from the adjacency of $S_{1,\infty}$ (resp.,
$T_\infty$) to $S_{1,n}$ (resp., $T_n$) for every $n\in\N$. 

We prove now that the singularities of classes  
$S_{n\geq 3}$ are non simple.
Set $\deg(\xi)=4$, $\deg(t)=1$.
The germ $f_a(\xi,t):=(\xi,t^6+t^2\xi+t^7+at^9)$, $a\in\R$, has 
a singularity of class $S_3$, due to Lemma \ref{lemma:2}.
By a direct computation, one verifies that
the weighted $9$-jets $f_a$ and $f_b$ are not $\A$-equivalent if $a\not=b$. 
Thus, in the neighboring of each such a map germ, there are 
infinitely many $\A$-orbits; that is, the germ is non simple.  
Let us consider a tangential family $f$ of class $S_3$. 
By Lemma \ref{lemma:2}, its prenormal form is 
$f(\xi,t)=(\xi,t^6+t^2\xi+\tilde \delta(6))$. 
It is easy to see that the weighted $9$-jet of $f$ is $\A$-equivalent to 
$(\xi,t^6+t^2\xi+A t^7+Bt^9)$ for some $A,B\in\R$. 
If $A\not=0$, we normalize it to $1$. 
Now the weighted $9$-jet of $f$ is of the form $f_a$ for some $a$, 
so $f$ is non simple. 
If $A=0$, we consider the deformation  
$F_\epsilon(\xi,t):=f(\xi,t)+(0,\epsilon t^7)$. 
For every $\epsilon\not=0$, $F_\epsilon$ is non simple, so $f$ is also
non simple. 
Finally, the singularities belonging to $S_{n>3}$ and $S_\infty$ are adjacent
to $S_3$, so they are not simple.  

It remains to prove that the singularities of class $U$ are non simple.
Let us set $f_a(\xi,t):=(\xi,t^4/4+2a
t^3\xi/3+t^2\xi^2/2)$.
Due to Lemma \ref{lemma:2}, $f_a$ is of type $U$ for every $a$.    
The vertical direction at the origin and the three branches of
the complex criminant set define four points on the projective complex
line, whose complex cross ratio is   
$2a^2-1+2 \sqrt{a^2-1}$.
This number is an invariant of tangential family germs, 
i.e. $\A$-equivalent tangential families have
the same cross ratio. 
Therefore, two germs whose $4$-jets are $f_a$ and $f_b$, 
are not $\A$-equivalent whenever their cross ratios are different,
i.e. $|a|\not=|b|$. 
Consider now the prenormal form $f(\xi,t)=(\xi,t^3\cdot \delta(1))$ of
any $U$-type tangential family. 
Set $F_\epsilon(\xi,t):= f(\xi,t)+(0,\epsilon(t^4+t^2\xi^2))$. 
For every $\epsilon$ small enough, $F_\epsilon$ is $\A$-equivalent to
a non simple map germ, so $f$ is non simple. 
This ends the proof of Theorem \ref{thm:1}.

\bigskip
\noindent {\bf Acknowledgments.}
The results of this paper are contained in my PhD Thesis \cite{these},
written under the supervision of V. I. Arnold, who proposed to me this
problem and careful read this paper.  
I also like to thank F. Aicardi, 
M. Garay and V. V. Goryunov for useful discussions and comments.


\bigskip

\noindent {\sc Gianmarco Capitanio} \\
\noindent {\it    Universit{\' e} D. Diderot -- Paris VII} \\ 
\noindent {\it    UFR de Math{\' e}matiques} \\
\noindent {\it    Equipe de G{\' e}om{\' e}trie et Dynamique}  \\
\noindent {\it    Case 7012 -- 2, place Jussieu} \\
\noindent {\it    75251 Paris Cedex 05} \\
{\it  e-mail:} Gianmarco.Capitanio@math.jussieu.fr

\end{document}